\documentclass[12pt]{amsart}
\usepackage{amssymb}
\usepackage{times,mathptmx}
\renewcommand{\baselinestretch}{1.5}

\begin{document}

\hfill Dedicated to V.S.Korolyuk on occasion of his  80-th
birthday

\title[]{Singular probability distributions \\ and fractal properties of sets of real
numbers \\ defined by the asymptotic frequencies of their s-adic
digits}

\author[S.Albeverio,  M.Pratsiovytyi, G.Torbin] {Sergio Albeverio$^{1,2,3,4}$,~  Mykola
Pratsiovytyi$^{5}$,~  Grygoriy Torbin$^{6,7}$}

\date{}

\begin{abstract}
    Properties of the set $T_s$ of "particularly  non-normal numbers" of the unit
    interval are studied in details ($T_s$ consists of real
    numbers $x$,  some of whose  s-adic digits have the
 asymptotic frequencies in the nonterminating $s-$ adic expansion
of $x$, and some do not). It is proven that the set $T_s$ is
residual in the topological sense (i.e., it is of the first Baire
category) and it is generic  in the sense of fractal geometry
($T_s$ is a superfractal set, i.e., its Hausdorff-Besicovitch
dimension   is equal to~1). A topological and fractal
classification of sets of real numbers via analysis of asymptotic
frequencies of digits in their s-adic expansions is presented.

\end{abstract}

\maketitle

{\renewcommand{\baselinestretch}{1}  $^1$ Institut f\"{u}r
Angewandte Mathematik, Universit\"{a}t Bonn, Wegelerstr. 6,
D-53115 Bonn (Germany) $^2$ SFB 611,  Bonn, BiBoS, Bielefeld -
Bonn $^3$ CERFIM, Locarno and Acc. Arch.,USI (Switzerland) $^4$
IZKS, Bonn; e-mail: albeverio@uni-bonn.de

$^{5}$ National Pedagogical University, Pyrogova str. 9, 01030
Kyiv (Ukraine); e-mail: prats@ukrpost.net

$^{6}$ National Pedagogical University, Pyrogova str. 9, 01030
Kyiv (Ukraine), $^{7}$ Institute for  Mathematics of NASU,
Tereshchenkivs'ka str., 3, 01601 Kyiv (Ukraine); e-mail:
torbin@imath.kiev.ua}



\section{Introduction}

Let us consider the classical $s-$adic expansion of $x\in[0,1]:$
$$ x = \sum\limits^\infty_{n=1} s^{-n}\alpha_n(x) = \Delta^s
\alpha_1(x) \alpha_2(x) ... \alpha_k(x)... ~,~~~ \alpha_k(x)\in
A=\{0,1, ... , (s-1)\}, $$
 and let
$N_i(x,k)$ be the number of digits "$i$" among the first $k $ digits
of the $s-$adic expansion of $x$, $i\in A.$ If the limit
$\nu_i(x)=\lim\limits_{k\rightarrow\infty} \frac {N_i(x,k)}{k}$
exists, then the number $\nu_i(x)$ is said to be the frequency of
the digit $"i"$ (or the asymptotic frequency of "$i$") in the
$s-$adic expansion of $x$.

A property of an element $x\in M$ is usually said to be "normal"
if "almost all" elements of $M$ have this property. There exist
many mathematical notions (e.g., cardinality, measure,
Hausdorff-Besicovitch dimension, Baire category) allowing us to
interpret the words "almost all" in a rigorous mathematical sense.
"Normal" properties of real numbers are deeply connected with the
asymptotic frequencies of their digits in some systems of
representation.

The set
$$N_s=\{x| (\forall i\in
A)~~\lim\limits_{k\rightarrow\infty} \frac {N_i(x,k)}{k} =
\frac{1}{s} \}
$$ is said to be the set of \textit{s-normal} numbers
(or the set of real numbers which are normal with respect to the
base s). It is well known (E.Borel, 1909), that the sets $N_s$ and
the set $N^*=\bigcap\limits_{s=2}^{\infty}N_s$ are of full
Lebesgue measure (i.e., they have Lebesgue measure 1).


The unit interval $[0,1]$ can be decomposed in the following way:
$$ [0,1]= E_s~\bigcup~D_s,$$
where
$$ E_s=\{x| \nu_i(x) ~~~exists, \forall i\in A  \},$$
$$
D_s=\{x| \exists i\in A , ~~\lim\limits_{k\rightarrow\infty} \frac
{N_i(x,k)}{k} ~~~does~~ not~~ exist \}.
$$
The set $D_s$ is said to be the set of non-normal real numbers.
Each of the subsets $E_s$ and  $D_s$ can be decomposed in the
following natural way.

The set
$$W_s=\{x| (\forall i\in A):~~\nu_i(x) ~~~exists, and ~~(\exists
j\in A): \nu_j(x)\neq \frac {1}{s}   \} $$ is said to be the set of
\textit{quasinormal} numbers. It is evident that $$E_s= W_s \bigcup
N_s,~~~ W_s \bigcap N_s=\varnothing.$$

The set $$L_s=\{x| (\forall i\in
A)~~\lim\limits_{k\rightarrow\infty} \frac {N_i(x,k)}{k} ~~~does~~
not~~ exist \}
$$ is said to be the set of \textit{essentially non-normal} numbers.

The set $$T_s = \{x| (\exists i\in
A):~~\lim\limits_{k\rightarrow\infty} \frac {N_i(x,k)}{k}
~~~does~~ not~~ exist,~ and ~ (\exists j\in
A):~~\lim\limits_{k\rightarrow\infty} \frac {N_j(x,k)}{k}
~~~exists \}$$ is said to be the set of \textit{particularly
non-normal} numbers.

It is evident that $$D_s= L_s \bigcup T_s, ~~~L_s \bigcap
T_s=\varnothing.$$

 The sets $N_s, W_s, T_s, L_s$ are everywhere dense sets, because
  the frequencies $\nu_i(x)$ do not depend on any finite
 number of s-adic symbols of $x$. It is also not hard to prove that
 these sets have the  cardinality of the continuum.

 The main purpose of the paper is to fill in completely the
 following table, which reflects the metric, topological and
 fractal properties of the corresponding sets:

 \begin{tabular}{|l|l|l|l|l|}
\hline
&Lebesgue measure&Hausdorff dimension & Baire category\\
\hline
$N_s$&&&\\
\hline
$W_s$&&&\\
\hline
$L_s$&&&\\
\hline
$T_s$&&&\\
\hline
\end{tabular}

 Let $\nu=(\nu_0, \nu_1, ... ,
\nu_{s-1})$ be a stochastic vector and let
$$
W_s[\nu]= \{ x:  x=\Delta^s \alpha_1(x) \alpha_2(x) ...
\alpha_k(x)...~,~~\lim\limits_{k \to \infty } \frac
{N_{i}^*(x,k)}{k}
 = \nu_i, \forall i \in A \}.
$$
The well known Besicovitch-Eggleston's theorem (see, e.g.,
\cite{Be,Eg}) gives the following formulae for the determination
of the Hausdorff-Besicovitch dimension $\alpha_0
\left(W_s[\nu]\right)$ of the set $W_s[\nu]$:
$$
 \alpha_0 \left(W_s[\nu]\right)= \frac { \sum\limits_{k=0 }^{s-1}
\nu_i \log \nu_i   } {- \log s}.
$$
From the latter formulae it easily follows that the set $W_s$ of all
quasinormal numbers is a superfractal set, i.e., $W_s$ is a set of
zero Lebesgue measure with full Hausdorff-Besicovitch dimension
($\alpha_0 \left(W_s\right)=1$).

 Properties of subsets of the set of non-normal
numbers have been intensively studied during recent  years (see,
e.g., \cite{PT,Ols1,Ols2,Ols3} and references therein). Some
interesting subsets of $D_s$ were studied in \cite{Ols1} by using
the techniques and results from the theory of multifractal
divergence points. In \cite{PT} it has been proven that the set
$D_s$ is superfractal.

In the paper \cite{APT3} of the authors it has been proven that the
set $L_s$ of essentially non-normal numbers is also superfractal and
it is of the second Baire category. Moreover, it has been proven
that the set $L_s$ contains an everywhere dense  $G_{\delta}$-set.
So, the sets $N_s, W_s, T_s$ are of the first Baire category.
 From these results it follows that  essentially non-normal numbers
 are generic in the topological
sense as well as in the sense of fractal geometry; nevertheless,
the set $L_s$ is small from  the  point of view of Lebesgue
measure.

The main goal of the present paper is the investigation of fractal
properties of the set $T_s$ of particularly non-normal numbers. To
this end we apply a probabilistic approach for the calculation of
the Hausdorff dimension of subsets. More precisely, we apply the
results of fine fractal analysis of singular continuous
probability distributions.

The first step of the fractal analysis of  a singular continuous
measure $\nu$ is the investigation of  metric, topological and
fractal properties of the corresponding topological support
 $S_{\nu}$ (i.e. the minimal closed set supporting the measure).
 These are good characteristics  only  for the class
  of  uniform Cantor-type singular measures. But, in general,
   they are  only "external
     characteristics", because  there exist essentially different singular
   continuous measures concentrating on the common topological support.
   The main idea of the paper \cite{APT3} consisted   in the construction of
   singular continuous measures whose topological supports coincide
   with some subsets of the set of essentially non-normal numbers.

  The second step of the fractal analysis of  a singular continuous
measure $\nu$ is the determination of the Hausdorff dimension
$\alpha_0
  (\nu)$ (and the local Hausdorff dimension ) of the measure, i.e., roughly speaking,
  finding  the Hausdorff dimension of the
  minimal (in the fractal dimension sense) supports (which are not necessarily
  closed) of the measure. This problem is much more complicated
  than the previous one (see, e.g., \cite{AT2}), especially in the
  case of essentially superfractal measures.

  In Section 2 we prove that for all $s\geq3$ the set $T_s$ is of full
   Hausdorff dimension.
   To prove the main result we  construct a sequence of singular continuous measures
   $\mu_p$ such    that the corresponding minimal dimensional
    supports consist of only particularly non-normal numbers, and
    apply the results of \cite{AT2} to perform a  fine fractal analysis of
    these supports.


\section{Fractal properties of the set of particularly non-normal numbers}

Let us study the sets $T_s$  of particularly non-normal numbers
which were defined in Section 1. It is easy to see that the set
$T_2$ is empty, because from the existence of the asymptotic
frequency $\nu_i(x)$ for some $i \in\{0,1\}$ the existence of
another  asymptotic frequency follows.

\textbf{Theorem 1.} \textit{For any positive integer $s\geq3$ the
set $T_s$ of particularly non-normal real numbers is superfractal,
i.e., the Hausdorff-Besicovitch dimension of the set $T_s$ equals
1.}

\proof
 To prove the theorem we shall construct a superfractal set $G\subset T_s$.

 In the sequel we usually shall not use the indices $s$ in the
notation of the corresponding subsets, since  $s$ will be an
arbitrary fixed  natural number greater than 2. Let us consider
the classical $s-$adic expansion of $x\in[0,1]:$ $x =
\sum\limits^\infty_{n=1} s^{-n}\alpha_n(x) = \Delta^s \alpha_1(x)
\alpha_2(x) ... \alpha_k(x)... ~.$ If $x$ is an s-adic rational
number, then we shall use the representation without the period
"$s-1$".

For a given $p\in N$ and for any $x\in[0,1)$ we define the
following mapping $\varphi_p$:

$$
 \varphi_p(x)= \varphi_p\left(\Delta^s \alpha_1(x) \alpha_2(x) ...
\alpha_k(x)...\right)= $$

$$=~ \Delta^s ~~ {\overbrace{00...0}^{s-1}}
{\overbrace{11...1}^{s-1}} ...
{\overbrace{(s-2)(s-2)...(s-2)}^{s-1}} (s-1)  \alpha_1(x)
\alpha_2(x)... \alpha_{s^2p}(x)
$$
$$
{\overbrace{00...0}^{2(s-1)}} ~~~{\overbrace{11...1}^{2(s-1)}} ...
{\overbrace{(s-2)(s-2)...(s-2)}^{2(s-1)}} (s-1)(s-1)
\alpha_{s^2p+1}(x) \alpha_{s^2p+2}(x)...
\alpha_{s^2p+2s^2p}(x)\ldots
$$
$$
{\overbrace{00...0}^{2^{k-1}(s-1)}}~~~~~~~~~~
 11...1
~...~~~ {\overbrace{(s-2)(s-2)...(s-2)}^{2^{k-1}(s-1)}}
{\overbrace{(s-1)(s-1)...(s-1)}^{2^{k-1}}}
\alpha_{(2^{k-1}-1)s^2p+1}(x) ...
\alpha_{(2^k-1)s^2p}(x)~~\ldots~~~.
$$
Let us explain the construction of $\varphi_p$. First of all we
divide the s-adic expansion of $x$ into groups in the following
way: the k-th group consists of the sequence
$(\alpha_{(2^{k-1}-1)s^2p+1}(x) ... \alpha_{(2^k-1)s^2p}(x)~~),
k\in N$. The s-adic expansion of $y=\varphi_p(x)$ is constructed
from the s-adic expansion of $x$ via inserting (before the k-th
group) the following series of fixed symbols $(0...01...1~~
...~~(s-2)...(s-2)~~(s-1)...(s-1) )$, where each symbol $"i"
(0\leq i\leq s-2)$ occurs $2^{k-1}(s-1)$ times, but the symbol
$"s-1"$ occurs $2^{k-1}$ times.

Let $M_p = \varphi_p([0,1))= \{y: y= \varphi_p(x), x\in [0,1)\}.$

For a given $p\in N$  and for any $y\in M_p$ we define the mapping
$\psi_p(y)$ in the following way: if $y=\varphi_p(x)= $
$$= ~ \Delta^s ~~ {\overbrace{00...0}^{s-1}}
 ... {\overbrace{(s-2)...(s-2)}^{s-1}}
(s-1)  \alpha_1(x) \alpha_2(x)... \alpha_{s^2p}(x)
$$
$$
{\overbrace{0...0}^{2(s-1)}} ~~~
...
{\overbrace{(s-2)...(s-2)}^{2(s-1)}} (s-1)(s-1) \alpha_{s^2p+1}(x)
\alpha_{s^2p+2}(x)... \alpha_{s^2p+2s^2p}(x)~~~~~~\ldots,
$$

then $z=\psi_p(y)=$
$$= ~ \Delta^s ~~ {\overbrace{0...0}^{s-1}}(s-1)
~~ ... {\overbrace{(s-2)...(s-2)}^{s-1}}(s-1)~~~~~~
(s-1)\left(01...(s-2)\right) \alpha_1(x) \alpha_2(x)...
\alpha_{s^2p}(x)
$$
$$
{\overbrace{00...0}^{(s-1)}}(s-1)
{\overbrace{00...0}^{(s-1)}}(s-1)
... {\overbrace{(s-2)(s-2)...(s-2)}^{(s-1)}} (s-1)
{\overbrace{(s-2)(s-2)...(s-2)}^{(s-1)}}(s-1)
$$

$$
 ~~~~ ~~~ (s-1)\left(01...(s-2)\right) ~~~
(s-1)\left(01...(s-2)\right)
 \alpha_{s^2p+1}(x)
\alpha_{s^2p+2}(x)... \alpha_{s^2p+2s^2p}(x)~~~~~~\ldots,~~ x\in
[0,1),
$$

i.e., the s-adic expansion of $z=\psi_p(y)$can be obtained from
the s-adic expansion of $y=\varphi(x)$ by using the following
algorithm:

1) after any fixed symbol "$(s-1)$" we  insert the following
series of symbols: $(01...(s-2))$;

 2) after any subseries
consisting of (s-1) fixed symbols $"i", 0\leq i \leq s-2$ we
insert the symbol $"s-1"$.

Let $f_p=\psi_p(\varphi_p)$ and let
$$
S_p=f_p([0,1))=\{z: z=f_p(x), x\in [0,1)  \} = \{z: z=\psi_p(y),
y\in M_p  \},
$$

$$
G_p=f_p([0,1))=\{z: z=f_p(x), x\in N_s  \}.
$$

The following two lemmas will describe some properties of the
constructed sets $G_p$.


{\bf Lemma 1.} \textit {For any $z= \sum\limits^\infty_{n=1}
s^{-n}\alpha_n(z)~~~\in G_p$ the limit
$\lim\limits_{n\rightarrow\infty} \frac {N_i(z,~~ n)}{n}$ does not
exist for any $i\in \{0,1, ... , s-2\}$, and
$\lim\limits_{n\rightarrow\infty} \frac {N_{s-1}(z,
~~n)}{n}~=~\frac{1}{s}.$}

\proof The set $G_p$ has the following structure:
$$
G_p=\{z: ~~z= ~ \Delta^s ~~ {\underbrace {\scriptstyle
{\overbrace{0...0}^{s-1}}(s-1) ~~ ...~~~
{\overbrace{(s-2)...(s-2)}^{s-1}}(s-1)~~ (s-1)(01...(s-2))
\alpha_1(x) \alpha_2(x)... \alpha_{s^2p}(x)}_{ first \,\,\,
group}}
$$

$$ {\underbrace {\scriptstyle
0...0(s-1)~~0...0(s-1) ... (s-2)...(s-2)(s-1)~ (s-2)...(s-2)(s-1)
~~(s-1)(01...(s-2)) ~ (s-1)(01...(s-2)) \alpha_{s^2p+1}(x)
\alpha_{s^2p+2}(x)... \alpha_{s^2p+2s^2p}(x)}_{second\,\,\,
group}}~
$$
$$~~\ldots~~~, x\in N_s~\}.$$ From $x\in N_s$ it follows that the symbol $"s-1"$
has the asymptotic frequency $\frac{1}{s}$ in the sequence
$\{\alpha_k(x)\}$ and the equality $\lim\limits_{n\rightarrow\infty}
\frac {N_{s-1}(z, n)}{n}~=~\frac{1}{s}$ follows from the
construction of the set $G_p$.

Let $l_k$ be the number of the position at which the above k-th
group of symbols ended, i.e., $l_k = s^2(p+1)(2^k-1).$

 Let $m_k'(i)$ be the number of the position at which the k-th
 series of the fixed symbols $"i"$ and $"(s-1)" (0\leq i\leq s-2)$
 ended, i.e., $m_{k+1}'(i)= s^2(p+1)(2^k-1) + s(i+1)2^k.$

 Let  $m_k''(i)$ be the number of the position at which the k-th
 series of the fixed symbols $"i" (0\leq i\leq s-2)$ started,
 i.e., $m_{k+1}''(i)= s^2(p+1)(2^k-1) + si2^k + 1$.

If $z\in G_p$, then there are $s(2^{k+1}-1)+ d_k)$ symbols $"i"
(0\leq i\leq s-2)$ among the first $m_{k+1}'(i)$ symbols of the
s-adic expansion of $z$, where $d_k$ is the quantity of the symbol
$"i"$ among the first $(2^k-1)s^2p $ s-adic symbols $\alpha_i(x)$ in
the expansion of $x=f_p^{-1}(z).$ Since $x$ is an s-normal number,
we have: $d_k= (2^k-1)sp+ \mathrm{o}(2^k).$

So, $\lim\limits_{n\rightarrow\infty} \frac {N_{i}(z,~~
m_{k+1}'(i))}{m_{k+1}'(i)}~=\lim\limits_{n\rightarrow\infty} \frac
{(2^{k+1}-1)s + (2^{k}-1)sp +s^{-1}
\mathrm{o}(2^k)}{s^2(p+1)(2^{k}-1)+s(i+1)2^k}~=~~
~\frac{p+2}{s(p+1)+i+1}.$

If $z\in G_p$, then there are $s(2^k-1)+ d_k$ symbols $"i"
~~(0\leq i\leq s-2)$ among the first $m_{k+1}''(i)-1$ symbols of
the s-adic expansion of $z$.

 So, $\lim\limits_{n\rightarrow\infty}
\frac {N_{i}(z,~~
m_{k+1}''(i)-1)}{m_{k+1}''(i)-1}~=~\lim\limits_{n\rightarrow\infty}
\frac {(2^{k}-1)s + (2^{k}-1)sp +s^{-1}
\mathrm{o}(2^k)}{s^2(p+1)(2^{k}-1)+si2^k}~=~~\frac{p+1}{s(p+1)+i}<
~\frac{p+2}{s(p+1)+i+1}.$

Therefore, for any $z\in G_p$ and for any $i\in\{0,1,..., s-2\}$
the limit $\lim\limits_{n\rightarrow\infty} \frac {N_i(z, n)}{n}$
does not exist.
\endproof

The following Corollary is immediate, using the definitions of
$G_p, T_s$ and Lemma 1:

 {\bf Corollary.} $G_p\subset T_s,\forall p\in N.$

\textbf{Lemma 2.} \textit{The Hausdorff-Besicovitch dimension of the
set $G_p$ is equal to $\frac {p}{p+2}$}.

\proof Let $B_p(i)$be the subset of $N$ with the following
property: $\forall k\in N, k\in B_p(i)  $ if and only if
$\alpha_k(f_p(x))=i$ for any $x\in [0,1)$, i.e., $ B_p(i)$
consists of the numbers of positions with the fixed symbols $"i"$
in the s-adic expansion of any $z\in S_p$.  ~~ Let
$B_p=\bigcup\limits_{i=0}^{s-1} B_p(i)$, and let $C_p=N\setminus
B_p$.

Let us consider the following random variable $\xi^{(p)}$ with
independent s-adic digits:
$$
\xi^{(p)} = \sum\limits_{k=1}^{\infty} s^{-k} ~\xi^{(p)}_k,
$$

where $\xi^{(p)}_k$ are independent random variables with the
following distributions: if $k\in B_p(i)$, then $\xi^{(p)}_k$
takes the value $"i"$ with probability 1. If $k\in C_p$, then
$\xi^{(p)}_k$ takes  the values $0,1,..., (s-1)$ with
probabilities $\frac{1}{s}, \frac{1}{s}, ..., \frac{1}{s}$.

It is evident that the set $S_p$ is the topological support of the
distribution of the random variable $\xi^{(p)}$. Actually, the
corresponding probability measure $\mu_p= P_{\xi^{(p)}}$ is the
image of Lebesgue measure on $[0,1)$ under the mapping $f_p =
\psi_p(\varphi_p)$, i.e., $\forall E\subset \mathcal{B}: \mu_p(E)=
\mu_p(E\bigcap S_p)= \lambda(f_p^{-1}(E\bigcap S_p))$.

\textbf{A)} Firstly we prove that $\alpha_0(G_p)\leq\frac {p}{p+2}.$
Since $G_p\subset S_p$, it is sufficient to show that
$\alpha_0(S_p)\leq\frac {p}{p+2}.$
 To this end we consider the
sequence $\{B_i^{(k)} \} ( k\in N, i\in \{1,2, ...,
s^{s^2p(2^{k-1}-1)}\})$ of special coverings of the set $S_p$ by
$s$-adic closed intervals of the rank
$m_k=l_k-2^{k-1}s^2p=s^2(p+1)(2^k-1)- 2^{k-1}s^2p $.
 For any $k\in N$ the covering $\{B_i^{(k)} \}$ consists of the
 $s^{s^2p(2^{k-1}-1)}$ closed s-adic intervals of $m_k$-th rank with
 length $\varepsilon_k = s^{-(s^2(p+1)(2^k-1)- 2^{k-1}s^2p)}.$

 The
$\alpha-$ volume of the covering $\{B_i^{(k)} \}$ is equal to
$$
l_{\varepsilon_k}^{\alpha}(S_p) = s^{s^2p(2^{k-1}-1)}\cdot
s^{-\alpha (s^2(p+1)(2^k-1)- 2^{k-1}s^2p)} =
s^{(p-\alpha(p+2))2^{k-1}s^2}\cdot~s^{\alpha(p+1)-p}.
$$

For the Hausdorff premeasure $h_{\varepsilon_k}^{\alpha}$ we have:
$h_{\varepsilon_k}^{\alpha}(S_p)\leq
l_{\varepsilon_k}^{\alpha}(S_p)$ for any $k\in N$. So, for the
Hausdorff measure  $H_{\alpha}$ we have $H_{\alpha}(S_p)\leq
\lim\limits_{k\rightarrow \infty} l_{\varepsilon_k}^{\alpha}(S_p) =
0$ if $\alpha > \frac {p}{p+2}.$

Hence, $\alpha_0(S_p)\leq\frac {p}{p+2}.$

\textbf{B)} Secondly  we prove that $ \alpha_0(G_p)\geq\frac
{p}{p+2}.$ To this end we shall analyze the internal fractal
properties of the singular continuous measure $\mu_p$.


For any probability measure $\nu$ one can introduce the notion of
the Hausdorff dimension of the measure in the following way:
 $$ \alpha_0(\nu)=
\inf\limits_{E \in N(\nu)} \{ \alpha_0(E), ~~ E \in \mathcal{B} \}
,$$ where $N(\nu)$ is the class of all "possible supports" of the
measure  $\nu$, i.e.,  $$ N(\nu)= \{ E: ~~ E \in \mathcal{B}, ~
\nu(E)=1 \} .$$

An explicit formula for the determination of the Hausdorff
dimension of the measures with independent Q*-symbols has been
found in \cite{AT2}. Applying this formula to our case ($q_{ik}=
\frac{1}{s}, \forall k\in N, \forall i\in \{0,1,..., s-1\}$), we
have
$$ \alpha_0(\mu_p)=\lim\limits_{\overline{ n \to
\infty}}~~~\frac{H_n}{n~lns}$$ where $H_n=\sum\limits_{j=1}^n h_j$,
and $ h_j$ are the entropies of the random variables $\xi_j^{(p)}:
~~h_j= - \sum\limits_{i=0}^{s-1} p_{ij} \ln p_{ij}. $

 If $j\in B_p$, then $h_j=0$. If $j\in C_p$, then $h_j=\ln s$.

 So,$$ \alpha_0(\mu_p)=\lim\limits_{\overline{ n \to
\infty}}~~~\frac{H_n}{n~lns}= \lim\limits_{\overline{ k \to
\infty}}~~~\frac{H_{m_k}}{m_k~lns}=  \lim\limits_{\overline{ k \to
\infty}}~~~\frac {s^2p(2^{k-1}-1)\ln s}{(s^2(p+1)(2^{k}-1)- p s^2
2^{k-1})~ \ln s}= \frac{p}{p+2}.$$

The above defined set $G_p=f_p(N_s)$ is a support of the measure
$\mu_p$, because $\mu_p = \lambda(f_p^{-1})$ and the Lebesgue
measure of the set $N_s$ of s-normal numbers of the unit interval
is equal to 1.

Since $G_p \in N(\mu_p)$ and $\alpha_0(\mu_p)= \frac{p}{p+2},$ we
get $ \alpha_0(\mu_p)\geq\frac{p}{p+2}$, which proves Lemma 2.
\endproof

{\bf Corollary.} \textit {The set $G_p$ is the minimal dimensional
support of the measure $\mu_p$, i.e.,$\alpha_0(G_p)\leq \alpha_0(E)$
for any other support $E$ of the measure $\mu_p$~.}

Finally, let us consider the set $G=\bigcup\limits_{p=1}^\infty
G_p.$  From Lemma 1 it follows that $G\subset T_s$. From Lemma 2 and
from the countable stability of the Hausdorff dimension it follows
that  $\alpha_0(G)= \sup\limits _p  \alpha_0(G_p) = 1.$ So,
$\alpha_0(T_s)=1$, which proves Theorem 1.
\endproof

Summarizing the results of Sections 1 and 2, we have for $s>2$:

\begin{tabular}{|c|c|c|c|c|} \hline
&Lebesgue measure&Hausdorff dimension & Baire category\\
\hline
$N_s$&1&1& first \\
\hline
$W_s$&0&1& first \\
\hline
$T_s$&0&1& first \\
\hline
$L_s$&0&1& second \\
\hline
\end{tabular}

\bigskip
For the case  $s=2$ we have a corresponding  result, but the
Hausdorff dimension of the set $T_s$ is equal to 0, because the
set $T_s$ is empty for $s=2$. \vskip 1 cm

\newpage
  \textbf{Acknowledgment}

This work was partly supported by DFG 436 UKR 113/78, DFG 436 UKR
113/80,   INTAS 00-257, SFB-611 projects and by Alexander von
Humboldt Foundation. The last two named authors gratefully
acknowledge the hospitality of the Institute of Applied Mathematics
and of the IZKS of the University of Bonn.

\end{document}